\documentclass[12pt]{article}
\usepackage{amssymb, url, amsthm, amsmath}
\newtheorem{theorem}{\noindent Theorem}

\newtheorem{lemma}{\noindent Lemma}
\newtheorem{definition}{\noindent Definition}

\newtheorem{proposition}{\noindent Proposition}

\newtheorem{question}{\noindent Question}
 \DeclareMathOperator{\Lip}{Lip}
\DeclareMathOperator{\const}{const} 
\DeclareMathOperator{\PLip}{PLip}

\author {J.~Melleray\thanks{%
University of Illinois at Urbana-Champaign. E-mail:{ \tt melleray@math.uiuc.edu}.},
F.~V.~Petrov\thanks{%
St.~Petersburg Department of Steklov Institute of Mathematics. E-mail:{ \tt fedorpetrov@mail.ru}.
Supported by the grant {NSh.4329.2006.1}},
A.~M.~Vershik\thanks{%
St.~Petersburg Department of Steklov Institute of Mathematics. E-mail:{ \tt vershik@pdmi.ras.ru}.
Supported by the grants {NSh.4329.2006.1,RFBR.05-01-00899}}}

\date{06.03.08}
\title{Linearly rigid metric spaces and the embedding problem.}

\begin{document}

 \maketitle
 \rightline{\it Dedicated to 95-th anniversary of L.V.Kantorovich}
\begin{abstract}
We consider the problem of isometric embedding of the metric spaces
to the Banach spaces; and introduce and study the remarkable class
of so called linearly rigid metric spaces: these are the spaces that
admit a unique, up to isometry, linearly dense isometric embedding
into a Banach space. The first nontrivial example of such a space
was given by R.~Holmes; he proved that the universal Urysohn space
has this property. We give a criterion of linear rigidity of a
metric space, which allows us to give a simple proof of the linear
rigidity of the Urysohn space and some other metric spaces. The
various properties of linearly rigid spaces and related
spaces are considered.
\end{abstract}

\section*{Introduction}
The goal of this paper is to describe the class of complete
separable metric (=Polish) spaces which have the following property:
there is a unique (up to isometry) isometric embedding of this
metric space $(X,\rho)$ in a Banach space such that the affine span
of the image of $X$ is dense (in which case we say that $X$ is
\textit{linearly dense}). We call such metric spaces
\textit{linearly rigid spaces}.

 The first nontrivial example of a linearly rigid space was Urysohn's
 universal space; this was proved by R.Holmes \cite{H}. Remember that
 the Urysohn space is the unique (up to isometry) Polish space which
 is universal (in the class of all Polish spaces) and
 ultra-homogeneous (= any isometry between finite subspaces extends
 to an isometry of the whole space). It was discovered by P.S.Urysohn
 in his last paper \cite{Ur}, which was published after his tragic
 death. The criterion of linear rigidity which we give in this paper
 is a weakening of well-known criterion of Urysohnness of the metric
 space, so the linear rigidity of Urysohn space is an evident
 corollary of our theorem.

 In this connection, in the first section we consider the general statement of the problem of
 isometric embedding of the metric spaces into Banach spaces.
 There exist several distinguish functorial embeddings of an arbitrary metric space
 into a Banach space. The first one is well-known Hausdorf-Kuratowski (HK):
 this is embedding into the Banach space of bounded continuous
 functions: $x \rightarrow \rho(x,.)$. Our point is to define isometric embedding using
 corresponding {\it norms and seminorms compatible
 with the metric} on the free space over metric space, or, in another
 words, - seminorms in the space $V(X)$ of all finite affine combinations of elements
 of $X$. This method was used in \cite{AE}, but it goes back to the idea of 40-th of
 the free group over metric or topological spaces.
 The most important compatible norm is the (KR)-norm, which based
 on the classical Kantorovich metric on the space of measures on the
 initial metric space; it was defined by L.V.Kantorovich in 1942 (\cite{K})
 in the framework of Kantorovich-Monge transportation problem.
 The Kantorovich-Rubinstein (KR) norm which is simply extension of Kantorovich metric to
 the space $V(X)$ (more exactly $_0(X)$ -see below) was defined in \cite{KR} for compact spaces.

 A remarkable fact is that (KR)-norm is the {\it maximal norm compatible with
 the given metric, so a linearly rigid space is a space for which
 (KR)-norm is the unique compatible norm}.

 The original definition of (KR)-norm is direct - as solution of
 transportation problem; the main observation by Kantorovich was that
 the dual definition of the norm used the conjugate space to the space $V_0$ which is
 the space of Lipschitz functions (with Lipschitz norm) on the metric space.
  Thus the completion of the space $V_0$ under that norm is predual
 to the Banach space $\Lip$ of Lipschitz function and is called sometimes
 as free Lipschitz space (see also \cite{DF,Kal,W}).  It also worth mentioning that geometry
 of the sphere in the (KR) norm is nontrivial even for
 finite metric space and is related to the {\it geometry of root
 polytopes of Lie algebras of series} $(A)$,
 and other interesting combinatorial and geometrical questions. This geometry
 is directly concerned with the problem of the embedding of the finite
 metric spaces to the Banach space; the authors do not know if this question had been discussed
 systematically in full generality.
 Because of maximality of KR-norm all other compatible norms and seminorms
 can be defined using some subspaces of the space of Lipschitz
 functions. We define the wider class of such norms article;
 the main example of this class is what we called the
 double-point - (dp)-norm; is used in the proof of the main theorem.

 In the second section we prove the main result ---  a criterion
 of linear rigidity in terms of distances. Namely, we prove that the
 characteristic property of that spaces is roughly speaking the
 following: any extremal (as the point of unit sphere) Lipschitz
 function of norm 1 is approximated by functions
 $\rho(x,\cdot)+\const$. The characteristic property of the Urysohn
 space is stronger: one does not need extremality of Lipschitz
 functions and there are some natural restrictions to the
 choice of constant. We prove in particular that the metric space is linearly
 rigid space {\it if only two norms -maximal (KR) and double-point (dp) are
 coincide.}
 We discuss some properties of the linearly rigid spaces, for example,
 we prove that such space must have infinite diameter; another property - the unit sphere
 of it is exclusively degenerated in a sense.

 In the third section we give several examples of linearly rigid
 spaces and first of all obtain the Holmes's result about linear
 rigidity of Urysohn space using apply our criterion and compare it
 with the criterion of Urysohnness. The integral and rational (with the
 distances more than 1) universal spaces are also linearly rigid. We
 discuss also the notions which are very close to rigidity. One of them
 is the notion of almost universal space which has an approximation property
 that is stronger than for linearly rigid spaces, but weaker
 than for the Urysohn space; Another one is the notion of weak rigidity of
 the metric spaces which corresponds to the the coincidence of $HK$ and $KR$
 norms which is weaker than linear rigidity.

 We formulate several questions appeared on the way. The geometry of
 those Banach spaces $E_X$ (and its unit spheres) which are
 corresponded to linearly rigid metric spaces $(X,\rho)$ are very
 intriguing. The most important concrete question: to define axiomatically the Banach space
 $E_{\Bbb U}$ which corresponds to Urysohn space $\Bbb U$, and can be called
 as Urysohn-Holmes-Kantorovich-Banach spaces or shortly - Urysohn Banach space.
 This is very interesting universal Banach space
 with a huge group of the linear isometries;
 it must be considered from various points of view.

\section{Isometric embeddings of metric spaces into Banach spaces.}

\subsection{Compatible norms and seminorms.}

 Let $(X,\rho)$ be a complete separable metric (=Polish)
 space\footnote{The main definitions, assertions and proofs of the
 paper are valid for nonseparable complete metric spaces}. Consider
 the free vector space $V={\Bbb R}(X)$ over the space $X$, and the
 free affine space $V_0={\Bbb R}_0(X)$ generated by the space $X$ (as
 a set) over the field of real numbers:

 $$
 V(X)={\Bbb R}(X)=\left\{\sum a_x\cdot \delta_x,\, x \in X, a_x \in {\Bbb R}\right\}
 $$
 $$\supset V_0(X)={\Bbb R}_0(X)=\left\{\sum a_x\cdot \delta_x,\, x \in X,
 a_x \in \Bbb R, \sum a_x=0\right\}
 $$
 (all sums are finite). The space $V_0(X)$ is a hyperplane in $V(X)$. We omit the mention of the space
 (and also of the metric, see below) if no ambiguity is possible; we will mostly consider only the
 space $V_0$. Another interpretation of the space ${\Bbb R}(X)$ (respectively, ${\Bbb R}_0(X)$) is
 that this is the {\it space of real measures} with finite support (respectively, the space of
 measures with finite support and total mass equal to zero: $\sum_x a_x=0$). Now we introduce the
 class of {\it seminorms} on $V_0$ {\it compatible with the metric $\rho$}. For brevity, denote by
 $e_{x,y}=\delta_x-\delta_y$ the elementary signed measure corresponding to an ordered pair $(x,y)$.

\begin{definition} {We say that a seminorm $\|\cdot \|$ on the space
$V_0(X)$ is compatible with the metric $\rho$ on the space $X$ if $\|e_{x,y}\|=\rho(x,y)$  for all
pairs $x,y \in X$.}
\end{definition}
Each compatible seminorm on $V_0$, can be extended to a seminorm on
the space $V$ by setting $\|\delta_x\|=0$ for some point $x \in X$,
and {\it vise versa}, the restriction of the compatible seminorm on
$V$ is a compatible seminorm on $V_0$; it is more convenient to
consider compatible seminorms only on the space $V_0$.

 The rays $\{ce_{x,y}\in V_0, \,c>0\}$ going through elementary signed measures will be
called {\it fundamental rays} (the set of fundamental rays does not
depend on the metric). If the metric is fixed, then a norm
compatible with this metric determines a unique vector of unit norm
on each fundamental ray; let us call these vectors (elementary
signed measures) {\it fundamental vertices corresponding to a given
metric}. They are given by the formula ($x \ne y$):
$$
\frac{e_{x,y}}{\rho(x,y)} \equiv {\bar e}_{x,y}.
$$
Thus the set of seminorms compatible with a given metric  $\rho$ is
the set of seminorms for which the fundamental vertices
corresponding to this metric are of norm one.

 The following useful elementary lemma describes all possible metrics on a space $X$
 in the geometrical terms of $V_0(X)$ .
\begin{lemma}
Let $X$ be a set. Consider the linear space $V_0(X)$ and specify
some points $c(x,y)\cdot e_{x,y}$ on the fundamental rays ${\Bbb
R}_+\cdot e_{x,y}$, where the function $c(x,y)$ is defined for all
pairs $(x,y)$, $x\ne y$,  positive, and symmetric: $c(x,y)=c(y,x)$.
This set of points is the set of fundamental vertices of some metric
on $X$ if and only if no point lies in the relative interior of the
convex hull of a set consisting of finitely many other fundamental
vertices and the zero.
\end{lemma}

\begin{proof}
 Consider the function defined by the formulas $\rho(x,y)=
c(x,y)^{-1}$, $x\ne y$, and $\rho(x,x)=0$. Let us check that the triangle inequality for this
function is equivalent to the property of convex hulls mentioned in the lemma. At first, if the
triangle inequality does not hold, say $\rho (c,a)>\rho(c,b)+\rho(b,a)$, then
$$
{\bar e}_{c,a}=\frac{\rho(c,b)}{\rho(c,a)}\cdot {\bar e}_{c,b}+ \frac{\rho(b,a)}{\rho(c,a)}\cdot
{\bar e}_{b,a}.
$$
The sum of coefficients $\frac{\rho(c,b)}{\rho(c,a)}+\frac{\rho(b,a)}{\rho(c,a)}$ is less then 1,
hence ${\bar e}_{c,a}$ lies in a relative interior of the triangle with vertices ${\bar e}_{c,b}$,
${\bar e}_{b,a}$ and 0.

Next, assume that $\rho$ is a metric. We need to prove that the sum $\lambda=\sum \lambda_i$ of the
coefficients in the linear representation
$$
{\bar e}_{x,y}=\sum_{i=1}^{n-1} \lambda_i{\bar e}_{x_i,y_{i}},\,\lambda_i\ge 0
$$
is at least $1$. Integrate the function $\rho(\cdot,y)$ with respect to measures in both sides. Due
to triangle inequality,
$$
\int \rho(\cdot,y) d {\bar e}_{x_i,y_i}=\frac{\rho(x_i,y)-\rho(y_i,y)}{\rho(x_i,y_i)}\le 1,
$$ so we get $1\le \sum \lambda_i$.
\end{proof}

Let ${\hat V}_0$ and $\hat V$ be the quotients of the spaces  $V_0$
and $V$ over the kernel $K=\{v:\|v\|=0\}$ of the seminorm and let
$\bar V$ and ${\bar V}_0$ be the completions of the spaces ${\hat
V}_0$ and $\hat V$ with respect to the norm.

\begin{proposition}
Every compatible seminorm $\|\cdot \|$ on the space $V_0(X)$ defines
an isometric embedding of the space  $(X,\rho)$ into the Banach
space $({\bar V}, \|\cdot \|$. Every isometric embedding of the
metric space $(X,\rho)$ into a Banach space $E$ corresponds to a
compatible seminorm on $V_0(X)$.
\end{proposition}

Indeed, obviously, the metric space $(X,\rho)$ has a canonical
isometric embedding into $\bar V$, and, conversely, it is easy to
see that if there exists an isometric embedding of the space
$(X,\rho)$ into some Banach space $E$, then the space $V_0$ is also
imbedded to the space $E$ linearly and the restriction of the norm
onto $V_0$ defines a compatible seminorm (not norm in general!) on
$V_0$.
\begin{question}
Let us call metric space solid if each compatible seminorm is a
norm. What metric spaces are solid? For example, what finite metric
spaces are solid? A more concrete question: what is the minimal
dimension of the Banach space into which a given finite metric space
can be isometrically embedded?
\end{question}

A similar (but different) question is studied in \cite{Sh}.

Of course it is enough to consider the case when the closure of the
affine hull of the image of $X$ in $E$ or the closure of the image
of space $V$ is dense in $E$. We will say that in this case the
isometric embedding of $X$ into a Banach space $E$ is {\it linearly
dense}. Thus, our problem is {\it to characterize the metric spaces
for which there is only one compatible norm or, equivalently, there
is a unique, up to isometry linearly dense embedding into a Banach
space.}

\subsection{Examples, functorial embeddings}

We will start with several important examples of isometric
embeddings and compatible norms.
\subsubsection{Hausdorf-Kuratowski embedding}
 The following is a well-known isometric embedding of an arbitrary metric space into a Banach
space:

\begin{definition} Define a map from the metric space $(X, \rho)$ into the Banach space ${\bar C}(X)$
of all bounded continuous functions on $X$, endowed with the sup-norm, as follows:
    $$X \ni x \mapsto \rho(x,.), V_0 \ni e_{x,y}\mapsto \rho(x,.)-\rho(y,.).$$
 We call it the Hausdorf-Kuratowski (HK) embedding.
 \end{definition}

 It is evident that this embedding is an isometry. In general, it is not a linearly dense
 embedding, because the image of the space $V(X)$ consists of very special
 Lipschitz functions. It is difficult to describe exactly the closed subspace of ${\bar C}(X)$
 that is the closed linear hull of the image under this embedding. At the same time, the
 corresponding compatible norm is given explicitly: for $v=\sum_k c_k\delta_{x_k}, \sum_k c_k=0$;
 $$
 \left\|\sum_k c_k \delta_{x_k}\right\|=
 \sup_z \left|\sum_k c_k \rho(z,x_k)\right|
 $$

\subsubsection{Compatible norms which are defined by the class of Lipschitz functions.}
Let us give another example of a class of embeddings, which will
play an important role. Choose the class of Lipschitz functions
${\cal L}\subset \Lip(X)$ which contains for each pair of points $x,y
\in X$ the Lipschitz function $f_{x,y}(\cdot)$ such that
$f_{x,y}(y)-f_{x,y}(x)=\rho(x,y)$

\begin{definition}
For every element $v=\sum_k c_k \delta_{z_k} \in V_0$ set
$$N_{\cal L}(v) = \sup_{f \in \cal L}|\sum c_k f(z_k)|.$$
Then $N_{\cal L}$ is a seminorm on $V_0$. We call it the  ${\cal
L}$-seminorm.
\end{definition}

 By definition of the  $N_{\cal L}$ we have $N_{\cal L}(e_{x,y})=\rho(x,y)$,
 so this seminorm is compatible with the metric $\rho$. Note that any compatible
 seminorm can be obtained in a similar way: by choosing a subset of Lipschitz
 functions $\cal L$.

 Let us consider the following important specific example of seminorms
 (at fact, they are actually norms) which is defined
 in this way. Let $\cal L$
be the following class of 1-Lipschitz functions: $${\cal L}=\{
\phi_{x,y}(\cdot)=\frac{ \rho(y,\cdot)-\rho(x,\cdot)}{2}, x,y \in X;
x \ne y\}.$$ Thus, the corresponding compatible norm is the
following:
$$\|v\|=\sup_{x,y} |\sum_k c_k \phi_{x,y}(z_k)|,$$ or in the notation for
continuous measures $\mu$:
$$\|\mu\|=\sup_{x,y}|\int \phi_{x,y}(z)d\mu(z)|.$$

\begin{definition}
We will call this norm on $V_0(X)$ double-point norm and denote it
as $\|\mu\|_{dp}$.
\end{definition}
We will heavily use double-point norm in the proof of the main
result.

 The choice of functions $\phi_{x,y}$ in the definition above
of double-point norm can be extended, for example, as follows: one may define
$$
\phi_{x,y}=\theta(x,y)\cdot \rho(x,\cdot)-(1-\theta(x,y))\cdot \rho(y,\cdot)
$$
for any function $\theta: X^2\rightarrow (0,1)$ such that
$\inf_{x,y} \theta(x,y)>0$ and $\sup_{x,y} \theta(x,y)=1$.

\subsubsection{Kantorovich embedding and maximal compatible norm}
Now we consider the most important compatible norm: it appeared as
consequence of the classical notion of Kantorovich (transport)
metric (\cite{42}) on the set of Borel measures on the compact
metric spaces later called as Kantorovich-Rubinstein norm \cite{KR}.

 The shortest way to define it is the following:
\begin{definition}
The Kantorovich-Rubinstein norm is the norm of $N_{\cal L}$-type
where the set $N_{\cal L}$ is the set of all Lipschitz functions
with Lipschitz norm one; more directly: let $v=\sum_k c_k
\delta_{z_k} \in V_0$ then $$\|v\|=\sup_u|\sum_k u(z_k)c_k|,$$ where
$u$ run over all Lipschitz functions over norm 1:
 $$\|u\|=\sup_{x\neq y} \frac{|u(x)-u(y)|}{\rho(x,y)}.$$
 We can define this norm not only on $V_0$ but on the space of all
 Borel measures $\mu$ on the metric space $(X,\rho)$ with compact
 support:
    $$\|\mu\|=\sup_{\|u\|=1} |\int_X u(z)d\mu(z)|.$$
\end{definition}

Immediately from the definition and from the remark above we have
the following theorem.
\begin{theorem}
 Let $(X,\rho)$ is a Polish space; then the KR-norm is the maximal compatible
norm; we will denote it as $\|.\|_{max}$.

In the geometrical terms, this means that the unit ball in KR-norm
is closed convex hull of the set of fundamental vertices ${\bar
e}_{x,y}$ (see above); the unit ball with respect to every seminorm
compatible with the metric $\rho$ contains the unit ball in KR-norm.
\end{theorem}

It is not difficult to prove this theorem directly. Recall that
initial definition of Kantorovich-Rubinstein norm was different and
the previous theorem is the duality definition of the norm with
comparison to original, in that setting the definition above is the
duality theorem due to L.Kantorovich (\cite{42}). As we mentioned
the conjugate space to the $V_0$ with KR-norm is the Banach space of
Lipschitz functions.

Denote the Banach space that is the completion of the space $V_0(X)$
with respect to the KR-norm by $E_{X,\rho}$. Sometimes this space
called "free Lipschitz space" (\cite{DF}. It is easy to check that
{\it The correspondence
$$(X,\rho)\mapsto E_{X,\rho} $$ is a functor from the category of
metric spaces (with Lipschitz maps as morphisms) to the category of
Banach spaces (with linear bounded maps as morphisms).}

Recall the initial definition was

\begin{definition}
Let $(X,\rho)$ be a Polish space. For the Borel measure
$\mu=\mu_+-\mu_-$, where $\mu_+$ and $\mu_-$ are the Borel
probability measures the Kantorovich-Rubinstein norm is defined as

$$\|\mu\|=\inf_{\psi}\int_X\int_X \rho(x,y)d\psi(x,y)$$
where $\psi$ runs over set $\Psi=\Psi(\mu_+,\mu_-)$ of all Borel
probability measures on the product $X\times X$, with the marginal
projection onto the first (second) factor equal to $\mu_+ (\mu_-)$.
\end{definition}

The equality between two definitions of the KR-norm is precisely
duality theorem in linear programming.

\subsubsection{Comments}

1.We have seen that the norm from Definition 6 defines a metric on
the simplex of probability Borel measures: $\|\mu\|\equiv
k_{\rho}(\mu_+,\mu_-)$, and just this metric was initially defined
by L.V.Kantorovich \cite{K}.\footnote{There are many authors who
rediscovered this metric later; unfortunately some of them did not
mention the initial paper by L.V.Kantorovich (see the survey
\cite{Met}).} The same is true for all other compatible norms: each
of which defines the metric on the affine simplex
$$V^+_0(X)=\{v=\sum_k c_k\delta_{z_k}\in V_0(X): c_k>0,\sum_k c_k=1\}$$ of
probability measures by formula:
$k_{\|.\|}(v_1,v_2)=\|v_1-v_2\|,v_1,v_2 \in V^+_0$.

More generally we can define the notion of {\it compatible metric}
on $V^+_0(X)$:

\begin{definition} A metric on the simplex $V^+_0(X)$, that is
convex as a function on affine set $V^+_0(X) \times V^+_0(X)$, and
which has the property $k_{\rho}(\delta_x,\delta_y)=\rho(x,y)$
is called compatible with the metric space $(X,\rho)$.
\end{definition}

There are many compatible metrics which do not come as above from
the compatible norms. F.e. the $L_p$-analogs, $p>1$ of the
Kantorovich metric are not generated by compatible norm. The
compatible metrics are very popular now in the theory of
transportation problems - see \cite{Vil}.

2.Our examples of compatible norms e.g. Hausdorf-Kuratowski,
double-point norm, maximal (Kantorovich-Rubinstein) norms etc. are
functorial norms (with respect to isometries as morphisms)in the
natural sense.

3.In opposite to existence of the maximal compatible norm there is
no {\it the least compatible norm}; moreover it happened that for
given compatible norm the infimum of the norms which which are less
that given norm, is {\it seminorm}, but not norm.

4.The unit ball in KR norm for finite metric space with all
distances equal to 1 is nothing more than generalization of the root
polytopes of the Lie algebras of series A. More generally, unit ball
for the general metric also could be considered as a polytope in
Cartan subalgebra.

\section{Main theorem: criterion of linear rigidity}

\subsection{Criterion of Linear rigidity}
Let $(X,\rho)$ be a metric space. We denote the Banach space of the
individual Lipschitz functions on the metric space $(X,\rho)$ as
$\Lip(X, \rho)$ or $\Lip(X)$ and denote quotient Banach space
$\PLip(X)=\Lip(X)/\{\const\}$ with norm: $\sup_{x \ne y}
\frac{|u(x)-u(y)|}{\rho(x,y)}$. The image of function $u\in Lip(X)$
under projection $\pi: \Lip(X) \rightarrow \PLip(X)$ we denote as
$\hat u$.
\begin{definition}

1)The Lipschitz function of type $\phi_x(\cdot)=\rho(x,\cdot)$ for
some point $x \in X$ is called a distance function,

2)The Lipschitz function $u$ is called admissible if
$$|u(x)-u(y)|\leq \rho(x,y)\leq u(x)+u(y).$$

3) Let $F\subset X$ (in particular $F=X$) and $u(\cdot)$ be an
individual Lipschitz function on metric space $F$ (equipped with the
induced metric) of norm 1. We say that the function $u$ is
representable (respectively, $\epsilon$-representable for $\epsilon
>0$) in $X$ if there exists a point $x \in X$ such that $u(z)=\rho(x,z)$
(correspondingly $|u(z)-\rho(x,z)|<\epsilon$) for all $z\in F$.

4)We say that $u$ is additively representable
($\epsilon$-representable) if there exists a constant $a\ in \Bbb R$
such that the function $u(\cdot)-a$ is representable
($\epsilon$-representable).
\end{definition}

Now we formulate the criterion of linearly rigidity of a metric
space.

Recall that for finite metric space $F$ a 1-Lipschitz
function $f$ on $F$ is called \emph{extremal}, if it is
an extreme point of the unit ball of Lipschitz functions (factored
by the constants). In other words, if functions $f\pm g$
are both 1-Lipschitz, then $g$ must be a constant function.

\begin{theorem}\label{criterion} The following assertions are
equivalent:

1)The complete metric space $(X,\rho)$ is linearly rigid;  that is,
 all norms on the space $V_0(X)$ compatible with metric coincide;

2)$\|\cdot \|_{max}=\|\cdot \|_{dp}$; --- i.e. two norms,  the
Kantorovich-Rubinstein norm and double-point norm, (see item
1.2.2-3) coincide;

3)(criterion) For each finite subset $F \subset X$ (with the induced
metric) and each $\epsilon>0$, any \textbf{extremal} 1-Lipschitz
function $u$ on the space $F$ is additively
$\epsilon$-representable.

4) The weak$^*$ closure of the convex hull of the set of distance
functions is the unit ball in the space $\PLip(X)$.\footnote{The last
formulation was suggested by one of the reviewer of the paper.}
\end{theorem}

 Remember that weak$^*$-topology on the space is defined by duality
 between the spaces $V_0(X)$ and $\PLip(X)$.

\begin{proof}
1) $\Longrightarrow$ 2) is trivial;

2) $\Longrightarrow$ 3): Let $(X,\rho)$ be a metric space, for which
maximal and double-point norms coincide. Our goal is to prove that
for any given finite $F\subset X$ and any extremal function $f$ on
the set $F$ there exists a point $x\in X$ and a constant $a$ such
that $\sup_{y\in F} |\rho(x,y)-f(y)-a|<\epsilon$.

If $F$ contains 1 or 2 elements, we may choose $x$ equal to one of these elements, so
let $F=\{x_1,\,x_2,\,\dots,x_n\},\,n\ge 3$.

Let us define a directed graph on $F$ for $n>2$ as follows: draw an edge
 $x_i\rightarrow x_j$ if $f(x_i)-f(x_j)=\rho(x_i,x_j)$.
 Note that the constructed graph regarded as an undirected
graph is connected. Indeed, if it is not connected, then for some
disjoint nonempty sets $F_1,\,F_2$ such that $F=F_1\cup F_2$ there
are no edges between $x_i$ and $x_j$ for $x_i\in F_1,\,x_j\in F_2$
Then for small positive $\epsilon$ the functions $f\pm
\epsilon\chi_{F_1}$ are 1-Lip too, this contradicts extremality of
$f$. (Here $\chi_{F_1}$ takes value 1 on $F_1$ and 0 on $F_2$). Let
us define an element $\mu\in V_0(F)$ as the sum
$$\mu:=\sum e_{a,b}$$
(here and up to end of the proof the summation is taken by all edges $a\rightarrow b$
of $F$).
A linear functional $\nu\rightarrow \int f d\nu$ attains its supremum on each function
$e_{a,b}$ for any edge $a\rightarrow b$. Hence it attains its maximum also on the sum of these
measures, i.e. on $\mu$. So
 $$\|\mu\|_{max} =\int f d\mu=\sum \rho(a,b).$$

Then also
$$
\|\mu\|_{dp}=\sum \rho(a,b)
$$

Hence for any $\epsilon>0$ there exist $x,y\in X$ such that the function $\phi_{x,y}$ satisfies the
following inequality:
$$
\left|\sum \rho(a,b)-\sum
\left(\phi_{x,y}(a)-\phi_{x,y}(b)\right)\right|<\epsilon.
$$
We have
$$
\rho(a,b)-(\phi_{x,y}(a)-\phi_{x,y}(b))=\left(\rho(a,b)-\rho(x,a)+\rho(x,b)\right)/2
+\left(\rho(a,b)+\rho(y,a)-\rho(y,b)\right)/2.
$$
Both summands are nonnegative and so both are less than $\epsilon$ for any edge $a\rightarrow b$.
It means that the function $g(\cdot)=\rho(x,\cdot)-f(\cdot)$ is almost constant on $F$ (since
$g(a)-g(b)$ is small for any edge and the graph is connected). So, $x$ satisfies the necessary
conditions.

\smallskip
3)$\Longrightarrow$ 1). We prove that if 3) holds for any finite $F$, then for every signed measure
$\mu\in V_0(F)$ and every norm $\|\cdot\|$ on $V_0(X)$ compatible with the metric, one has
$\|\mu\|=\|\mu\|_K$.

Recall that the unit ball of $\|\cdot\|_K$-norm is the closed convex hull of the
points ${\bar e}_{a,b}$; every
finitely-supported measure $\mu$ such that $\|\mu\|_K=1$ is a convex combination of some ${\bar
e}_{a_k,b_k}$. Applying this to the measure $\mu/\|\mu\|_K$ we get
$$
\mu=\sum_{k=1}^N \alpha_k {\bar e}_{a_k,b_k}, \quad a_k,b_k\in F,\;\alpha_k\ge
0,\quad\|\mu\|_K=\sum \alpha_k.
$$

The points ${\bar e}_{a_k,b_k}$ lie on some face of the unit ball of the space $E_{F}$. We may
assume without loss of generality that it is a face of codimension~$1$. The corresponding
supporting plane is determined by some linear functional of norm $1$, i.e., some $1$-Lipschitz
function $f$ on $F$. Then for every $k$, $f({\bar e}_{a_k,b_k})=1$, that is $f(a_k)-f(b_k)=\rho(a_k,b_k)$.
$f$ is an extremal Lipschitz function on $F$.

First consider the case when all the $a_k$ are equal: $a_k=a$.

By assumption, there is  a point $c\in X$ such that $\rho(c,a)\ge
\rho(c,b_k)+\rho(a,b_k)-\varepsilon$.

We have
$$
\|\mu\|=\left\|\sum \alpha_k {(\delta_a-\delta_c)+(\delta_c-\delta_{b_k})\over
\rho(a,b_k)}\right\|$$
$$\ge \sum \alpha_k\cdot {\rho(a,c)\over \rho(a,b_k)}-\sum \alpha_k
\cdot {\rho(c,b_k)\over \rho(a,b_k)}\ge \sum \alpha_k\left(1-{\varepsilon\over \min_k
\rho(a,b_k)}\right).
$$

Letting $\varepsilon \to 0$, we obtain
$$
\|\mu\|=\sum \alpha_k=\|\mu\|_K.
$$
So, for any measure of the type $$\mu=\sum \alpha_k {\bar e_{a,b_k}},\,\alpha_k\ge 0$$ we get
$\|\mu\|=\|\mu\|_K=\sum \alpha_k$ (it is easy to check that $\|\mu\|_K=\sum \alpha_k$ for any such
measure).

Now consider the general case. Find a point $d$ such that $\rho(d,a_k)\ge
\rho(d,b_k)+\rho(a_k,b_k)-\varepsilon$. Then we obtain
\begin{eqnarray*}
\|\mu\|&=&\left\|\sum \alpha_k {(\delta_{a_k}-\delta_d)+(\delta_d-\delta_{b_k})\over
\rho(a_k,b_k)}\right\| \\
&\ge& \left\|\sum \alpha_k {\delta_{a_k}-\delta_d\over \rho(a_k,b_k)}\right\|-\sum \alpha_k
\left\|{\delta_d-\delta_{b_k}\over
\rho(a_k,b_k)}\right\| \\
&\ge&\sum \alpha_k {\rho(a_k,d)\over \rho(a_k,b_k)}-\sum \alpha_k {\rho(d,b_k)\over
\rho(a_k,b_k)}\ge \sum \alpha_k+o(1)
\end{eqnarray*}
(in the second inequality we use the case which we have considered at first). This completes the
proof in this general case too.

3) $\Longrightarrow$ 4) and 4) $\Longrightarrow$ 3).

 It suffices to prove that any {\it extremal} Lipschitz function
 lies in the weak$^*$-closure of
the set distance functions (which consists of extremal Lipschitz
functions itself). Take the Lipschitz function $u \in PLip(X)$ which
is an extremal on the unit ball of the space. It lies in the
weak$^*$-closure of a set $L$ iff any weak$^*$-open neighborhood
$W\ni u$ intersects $L$. But such a neighborhood by definition of
weak topology is defined by a finite subset of the space $X$, so we
just get the criterion in formulation 3).
\end{proof}
Remark that it is possible to give a shorter proof of $1)\Leftrightarrow 3)$ directly
by using the Krein-Milman theorem; but the proof above is elementary
and reduces the general case to the case of a finite metric space
and did not use infinite constructions. We will compare our
criterion with the characteristic property of the Urysohn space in
section 3.

\textbf {Remark.} We see that the coincidence of the double-point
norm with the Kantorovich norm implies coincidence of all compatible
norms. But instead of double-point norm we can choose any norm which
was defined in 1.3.1. The proof is essentially the same. So there
are many ways to define a norm, compatible with the metric, the
coincidence of which with Kantorovich norm implies the linear
rigidity.

\subsection{Properties of linearly rigid spaces}

\subsubsection{Unboundedness of linearly rigid spaces.}
The following theorem shows that a linearly rigid space cannot have
a finite diameter if it has more than two points.

\begin{theorem} A linearly rigid metric space
$X$ containing more than two points is of infinite diameter and, in particular, noncompact.
\end{theorem}

\begin{proof} Assume the contrary. Without loss of generality
we assume that the space $X$ is complete. Fix a point $a\in X$, denote by $r_a$ the supremum of the
distances $\rho(a,x)$ over $x\in X$, and choose a sequence of points $(x_n)$ such that
$\rho(a,x_n)\ge r_a-1/n$. Then pick a countable dense subset $\{y_n\}$ of $X$, and define a
sequence $(z_n)$ by setting $z_{2n}=x_n$ and $z_{2n+1}=y_n$. Consider the points
 ${\bar e}_{a,z_k}$, $k=1,{\ldots} ,N$. They lie on the same face
of the unit ball of the space $E_{X_N}$, where
$X_N=\{a,z_1,z_2,\dots,z_N\}$. Applying Theorem \ref{criterion}, we
may find for each $N$ a point $c_N$ such that
$$\rho(a,c_N)\ge \rho(a,z_k)+\rho(z_k,c_{N})-1/N,\qquad k\le N.$$

In particular, $\rho(a,c_N) \ge \rho(a,x_k)+\rho(x_k,c_{N})-1/N,\
2k\le N.$ Hence $\rho(a,c_N)\rightarrow r_a$ and therefore
$\rho(x_k,c_{N})\to 0$ as $k,N\rightarrow \infty$, so that the
sequences $(x_k)$, $(c_{k})$ are both Cauchy and have a common limit
$a'$. The point $a'$ satisfies the equalities
$\rho(a,x)+\rho(x,a')=\rho(a,a')$ for all $x\in X$ (this is why we
used the countable dense set $\{y_n\}$ in the definition of our
sequence $(z_n)$).

Such a construction may be done for any point $a \in X$; note that
for any $a,b\in X$ such that $a\ne b\ne a'$ (such $a,b$ do exist if
$X$ has more than two points) we have
$2\rho(a,a')=\rho(a,b)+\rho(a,b')+\rho(a',b)+\rho(a',b')=2\rho(b,b')$,
whence $\rho(a,a')\equiv D<{\infty}$. It also follows that
$\rho(a,b)=\rho(a',b')$.

Without loss of generality, $\rho(a,b')=\rho(b,a')\ge \rho(a,b)=\rho(a',b')=1$. Let
$A=\{a,b,a',b'\}$. Define a function $\varphi$ by the formulas $\varphi(a)=\varphi(a')=1$,
$\varphi(b)=\varphi(b')=0$.

Such a function is 1-Lipschitz on $\{a,b,a',b'\}$; the corresponding face contains the points ${\bar
e}_{a,b}$, ${\bar e}_{a',b'}$. Hence there exists a point $c$ such that
$$
\rho(c,a')\ge \rho(c,b')+1/2,\qquad \rho(c,a)\ge \rho(c,b)+1/2
$$

We have
$$
\rho(a,a')=\rho(a,c)+\rho(c,a')\ge \rho(c,b')+\rho(c,b)+1=D+1.
$$
The obtained contradiction proves the theorem.
\end{proof}

\subsubsection{How to construct inductively a linearly rigid space}

Let $(X,\rho)$ is linearly rigid metric space and $E_{X,\rho}$
corresponding Banach space. The properties of the unit sphere of the
space $E_{X,\rho}$ are very peculiar and can be used for the
recursive construction of the metric space and corresponding Banach
space. We give a draft of the inductive construction; it is based on
the following finite-dimensional

\begin{theorem}[Piercing theorem] Let $(Y_1,r_1)$ be an arbitrary finite metric space,
$\epsilon>0$, and $\Gamma$ be a face of the unit ball of the space
$E_{Y_1,r_1}$ (e.g. $V_0(Y_1)$ with compatible norm). Then the space
$(Y_1,r_1)$ can be isometrically embedded into a finite metric space
$(Y_2,r_2)$ so that there exists a face $\Delta$ of the unit ball of
the space $E_{Y_2,r_2}$ containing $\Gamma$ and two vectors ${\bar
e}_{z_1,z_2}$ and ${\bar e}_{u_1,u_2}$ such that the line segment
connecting them intersects the face $\Delta$ at an interior point.
\end{theorem}

The proof is direct. Note that if an interior point of a
face is of norm one, then all points of the face are also of norm
one. Enumerating the sequences of faces of the root polytopes
already constructed and ``piercing'' them by new line segments, we
obtain a sequence of finite metric spaces for which all faces of all
root polytopes are rigid; hence the completion of the constructed
countable space will be linearly rigid. Note that such a degeneracy
of the unit sphere is typical for universal constructions (cf. the
Poulsen simplex - \cite{L}). We expect the the geometry of the unite
sphere of the spaces $E_{X,\rho}$ for linearly rigid metric space
$(X,\rho)$ is very unusual and interesting.

\section{Examples of linearly rigid spaces and related problems}
Up to now we haven't provided any example of linearly rigid space. A
trivial example of linearly rigid spaces is one- and two- point
spaces. R.~Holmes in \cite{H} had discovered that universal Urysohn
space has this property: each isometric embedding of it to the
universal space $C([0,1])$ generates as a linear hull the isometric
Banach spaces, consequently, in our terminology this means that
Urysohn space is linearly rigid, and this was the first nontrivial
example of such a space. We will deduce this result as well as other
examples as easy consequence of our criterion.

\subsection{Criteria of Urysohnness and linear rigidity of the Urysohn space.}

In order to prove the linear rigidity of the Urysohn universal space
we recall its characterization. The following criteria is a
treatment of the different characterization from the original papers
of Urysohn \cite{Ur} and subsequent papers \cite{K,V,G}.
\begin{theorem} {A Polish space $(X,\rho)$ is
isometric to the universal Urysohn space if and only if for every
$\epsilon >0$,  any every finite subset $F\subset X$, every
admissible Lipschitz function $u$ on $F$ is $\epsilon$-representable
in the space X}. (see Definition 8.4)
\end{theorem}

From the point of view of functional analysis we can reformulate
this condition in the following words: the set of distance function
is weakly$^*$ dense in the unit ball of the space $Lip(X)$.

\begin{theorem}[R.~Holmes \cite{H}]
The Urysohn space $\Bbb U$ is linearly rigid.
\end{theorem}

\begin{proof}
It suffices to compare the assumptions of the criterion of
universality above and linear rigidity criterion from the previous
section: the assumptions of the latter require additive
$\epsilon$-representability of extremal Lipschitz functions, while
the universality criterion requires $\epsilon$-representability of
all positive Lipschitz functions.

In other words, accordingly to item 4) of Theorem 2, the linear
rigidity is equivalent to weak$^*$ density of convex hull of the set
distance function in unit ball of $PLip(X)$ while the universality
is equivalent to the density of the set of distance functions in the
in the unit ball of $Lip(X)$, which is a much stronger condition.
\end{proof}

It is natural to call the Banach space $E_{\Bbb U}$ (completion of
the space $V_0(\Bbb U)$ with respect to the unique compatible norm)
as "Urysohn Banach space"\footnote{or Urysohn-Holms-Kantorovich
Banach space; in the case of coincidences of it with Gurariy space
-see Question 2 below - also add name Gurariy.}. The geometry of
this space seems to be very interesting. First of all $E_{\Bbb U}$
as Banach space is universal, which means that each separable Banach
space can be linearly isometrically embedded to it. As Professor
V.Pestov pointed out this follows from a strong theorem of Godefroy
and Kalton \cite{GK}, which states that if some separable Banach
space $F$ has an isometric embedding into a Banach space $B$, then
it also has a linear isometric embedding into $B$.

However, $E_{{\Bbb U}, \rho}$ is not a homogeneous universal as
Banach space -it is known that there is no separable Banach space,
in which each linear isometry between any two finite dimensional
isometric linear subspaces can be extended to a global isometry of
the space (\cite{Li,Gu,L}). The {\it Gurariy space} (\cite{Gu}) has
$\epsilon$-version of this property.

\begin{question}
Is the Urysohn Banach space $E_{\Bbb U}$ isometrically isomorphic to
the Gurariy space?
\end{question}

\subsection{The further examples}

\subsubsection{Rational discrete universal metric spaces are linearly rigid}

  Let us discuss other examples of linearly rigid universal spaces.

\medskip\noindent
 Let us consider the countable metric space denoted by $\Bbb Q\Bbb U_{\geq 1}$. It is a universal
and ultra-homogeneous space in the class of countable metric spaces with rational distances not
smaller than one. Such a space can be constructed in exactly the same way as the Urysohn space.

\begin{theorem}
The space $\Bbb Q\Bbb U_{\geq1}$ is linearly rigid.
\end{theorem}

Indeed, the assumptions of the criterion of linear rigidity are obviously satisfied.

This example, as well as the next one, is of interest because it is an example of a discrete
countable linearly rigid space. Thus the corresponding Banach space $E_{\Bbb Q\Bbb U_{\geq 1}}$ has
a basis. It is not known whether the space $E_{\Bbb U}$ has a basis.
\begin{definition}
Let us call metric space $(X,\rho)$\emph{almost universal} if the
set of distance functions is weakly$^*$ dense in the unit
ball of the space $\PLip(X)$.
\end{definition}
This notion is in between universal Urysohn and linearly rigid
spaces but does not coincide with any of them. It is easy to prove
\begin{proposition}
The space $\Bbb Q\Bbb U_{\geq1}$ is almost
universal.
\end{proposition}
In the same time integral universal space is linearly rigid (see next
item), but not almost universal.
\begin{question}
To describe all almost universal metric spaces.
\end{question}

\subsubsection{Integral universal metric space is linearly rigid}

The following example is of special interest also for another
reason. Consider the space $\Bbb Z\Bbb U$, the universal and
ultra-homogeneous space in the class of metric spaces with integer
distances between points. Let us show that it is also linearly
rigid. For this, let us check the condition of the criterion of
linear rigidity. Fix a finite set $X_n$ in the space $X$ and
extremal ray $L=\{\lambda f\}; \lambda>0 $ of the set of Lipschitz
functions on $X_n$. Note that the differences of the coordinates of
every vector from the ray $L$ are integers; hence on this ray there
is a vector with integer coordinates, which is realized as the
vector of distances between some point $x\in X$ and the points of
set $X_n$.

Let us introduce a graph structure on this countable space by
assuming that pairs of points at distance one are neighbors. This
graph has remarkable properties: it is universal but not homogeneous
(as a graph), its group of isomorphisms coincides with the group of
isometries of this space regarded as a metric space. As follows from
\cite{Cam, CamVer}, there exists an isometry that acts transitively
on this space.
\begin{question}
To characterize this graph using universality of it as the metric
space.
\end{question}

\medskip\noindent

\subsection{Weakly linearly rigid spaces}

We have proved that coincidence of KR and DP-norms leads to linear
rigidity. Now let us compare another two compatible norms: the
(KR)-norm and (HK)-norm.

\begin{definition}
A metric space for which the maximal (or KR-) and HK-norms coincide is
called a weakly linear rigid space (WLR-space).
\end{definition}

If the space $(X,\rho)$ is linearly rigid, then all compatible norms
must coincide with the maximal norm of $\sum c_k \delta_{x_k}$,
which is
$$
\|\sum c_k\delta_{x_k}\|_{HK}=\|\sum
c_k\delta_{x_k}\|_{KR}=\sup_{F\in \Lip_1(X)} \sum c_k F(x_k),
$$
where the supremum is taken by all 1-Lipschitz functions on $X$ (or,
equivalently, on the set of $\{x_k\}$). The equality above means
that every extremal Lipschitz function is almost realized for a
function $\pm \rho(x,\cdot)$ for some $x\in X$. The criterion of
linear rigidity says that any extremal Lipschitz function on $X$ is
almost realized as a function $\rho(x,\cdot)$. The only difference
is the absence of $\pm$. Recall that if we take a supremum not by
the set of functions $\pm \rho(x,\cdot)$, but by the set of
functions $\frac12 (\rho(x,\cdot)-\rho(y,\cdot))$, then the
coincidence of the corresponding norm with the maximal (KR) norm
implies linear rigidity. But this difference is quite essential.
There are some spaces which are weakly linearly rigid but not
linearly rigid. We do not have any less or more complete description
of such spaces. So only some examples follow.

1. Any metric space on three points gives a WLR-space (but it is never linearly rigid, as any finite
metric space on more than two points, see theorem 3).

2. Let us define a family of WLR-spaces on four points. Let points be $A,\,B,\,C,\,D$; let $a,\,b,\,c$ be
arbitrary positive numbers and define:
$\rho(D,A)=a,\,\rho(D,B)=b,\,\rho(D,c)=c,\,\rho(A,B)=a+b,\,\rho(B,C)=b+c,\,\rho(A,C)=a+c$. It is
easy to see that any extremal Lipschitz function is realized as $\pm\rho(X,\cdot)$ for some $X\in
\{A,B,C,D\}$.

\begin{question}
Does there exist the finite WLR metric space with more than four points?
\end{question}

It looks more likely that there are no such spaces at least with sufficiently many points.

3.In the same time there are an infinite WLR-space which is not be
linearly rigid either. Consider the Urysohn space $\Bbb U$ and fix a
point $a\in \Bbb U$. Add a point $a'$ to the space $\Bbb U$ and
define the distances $\rho(a',x)=\rho(a,x)+1$ for any $x\in X$ (in
particular, $\rho(a,a')=1$). This new space $U'=\Bbb U\cup \{a'\}$
is WLR, but not linearly rigid.

\subsection{Extremality and the properties of the
Banach spaces $E_{X,\rho}$}

The set of all possible distance matrices (semimetrics) is a convex
weakly (i.e. in pointwise topology) closed cone (see \cite{V}).
 If the distance matrix of a Polish space $X$, which corresponds to some
 dense sequence in $X$, belongs to extremal
 ray of this cone, we say that $X$ is \emph{extremal}. This property
 does not depend on the choice of dense sequence, so the
 definition is correct. This notion is interesting even for finite
 metric spaces (see.\cite{A}).
 \begin{question}
To describe extremal finite metric spaces with $n$ points, to
estimate exact number of such spaces or asymptotics on $n$.
 \end{question}

Note that the universal real Urysohn space $\Bbb U$ {\it is
extremal} metric space (\cite{V}). It follows from the genericity of
$\Bbb U$ that the distance matrices of everywhere dense systems of
points of extremal metric spaces form an everywhere dense
$G_{\delta}$-set in the space of distance matrices. The integer
space $\Bbb Z\Bbb U$ is also extremal; the extremality of both
spaces follows from a result of Avis \cite{A}, which states that
every finite metric space with commensurable distances (i.e. such
that the ratio of any two distances is rational number) can be
embedded into a finite extremal metric space, and hence the
assumptions of the criterion of linear rigidity are satisfied. Using
the criterion of linear rigidity given above, and the procedure of
the previous section one can build an example of a not extremal
linearly rigid metric space.

\textbf{Acknowledgment.} The authors are grateful to V.Pestov for
useful discussions and to the reviewers for valuable suggestions.

\end{document}